\documentclass[a4paper]{article}
\usepackage[utf8]{inputenc}

\usepackage[margin=1in]{geometry}

\usepackage{amsmath, amssymb, tcolorbox}
\usepackage[hidelinks]{hyperref}

\usepackage{xcolor}
\usepackage[dvipsnames]{xcolor}
\definecolor{blauw}{RGB}{24,65,131}
\definecolor{oranje}{RGB}{232,78,15}

\usepackage{amsthm}
\usepackage{thmtools}
\usepackage[noabbrev]{cleveref}
\newtheorem{thm}{Theorem}

\newtheorem{crl}[thm]{Corollary}
\newtheorem{prop}[thm]{Proposition}

\theoremstyle{definition}

\newtheorem{ex}[thm]{Example}

\renewcommand{\phi}{\varphi}

\newcommand{\FF}{\mathbb F}

\renewcommand{\leq}{\leqslant}
\renewcommand{\geq}{\geqslant}

\newcommand{\sett}[2]{ \left\{ #1 \, \, || \, \, #2 \right \} }

\newcommand{\floor}[1]{\left \lfloor #1 \right \rfloor}

\newcommand{\mh}{\mathcal H}

\newcommand{\ms}{\mathcal S}

 \DeclareMathOperator{\rk}{rk}

\makeatletter
\let\@fnsymbol\@arabic  
\makeatother

\title{A note on the classification of classical distance-regular graphs of negative type and the non-existence of hemisystems}
\author{
 Sam Adriaensen\thanks{Vrije Universiteit Brussel, Pleinlaan 2, 1050 Elsene, Belgium. \{\href{mailto:sam.adriaensen@vub.be}{sam.adriaensen}, \href{mailto:jan.de.beule@vub.be}{jan.de.beule}, \href{mailto:sam.mattheus@vub.be}{sam.mattheus}\}@vub.be.}
 \and 
 Jan De Beule\footnotemark[1]
 \and 
 Jozefien D'haeseleer\thanks{Ghent University, Krijgslaan 297, Building S8, 9000 Gent, Belgium. \href{mailto:jozefien.dhaeseleer@ugent.be}{jozefien.dhaeseleer@ugent.be}.}
 \and 
 Sam Mattheus\footnotemark[1]
}
\date{}

\begin{document}

\maketitle

\begin{tcolorbox}[colback=red!5!white, colframe=red!75!black, title={\scshape Disclaimer}]
 There are problems with the correctness of proofs in the literature cited in this paper.
 The reader should therefore consider the main conclusion of this note as {\scshape not proven} until these issues are resolved.
 More specifically:

 \begin{enumerate}
  \item In \cite[Theorem 1.1]{Hiraki}, Hiraki proved a theorem in the language of near polygons.
  \item In \cite[Proposition 11.3]{VanDam:Koolen:Tanaka}, the authors rephrased the theorem of Hiraki in terms of $m$-bounded distance-regular graphs.
  In particular, this gives a sufficient condition for a distance-regular graph to be $m$-bounded.
  \item In \cite[Lemma 2.4]{Koolen:Zhi}, the authors cite \cite[Theorem 1.1]{Hiraki} and \cite[Proposition 11.3]{VanDam:Koolen:Tanaka} incorrectly.
  They state that the condition for $m$-boundedness is necessary and sufficient instead of only sufficient.
  In other words, they cite an ``$m$-bounded if condition C'' statement as an ``$m$-bounded if and only if condition C'' statement.
  \item In \cite[Lemma 2.5]{Tian:Etal}, the authors use the incorrect statement of \cite[Lemma 2.4]{Koolen:Zhi}, and the error in this statement is crucial for the proofs in their paper. 
 \end{enumerate}
 We can thus not be sure of the correctness of the results in \cite{Tian:Etal}, and therefore neither of the conclusions drawn in this paper. 
\end{tcolorbox}

\begin{abstract}
 The goal of this note is to connect some interesting results in the literature on algebraic graph theory and finite geometry.
 In 1999, Weng gave an almost complete classification of classical distance-regular graphs of negative type with diameter at least 4.
 He proved that these graphs are either dual polar graphs of Hermitian polar spaces, Hermitian forms graphs, or fall into a last category.
 It was recently proved by Yian et al.\ that the latter category does not exist when the diameter equals 3, which by Weng's results proves that they do not exist for bigger diameter.
 Using a result of Vanhove, this proves that certain hemisystems in Hermitian polar spaces cannot exist.
\end{abstract}

The goal of this note is to combine the recent findings of Tian et al.\ \cite{Tian:Etal} with the existing literature, to draw some interesting conclusions.
We emphasize this note does not contribute any new mathematical ideas, rather we connect some dots in the realm of mathematical literature.

Let us start by presenting the main findings, the necessary mathematical background can be found in \Cref{Sec:Prel}.

Weng \cite{Weng} gave an almost full classification of classical distance-regular graphs of negative type.
We will use the common notation $a_i, b_i, c_i$ for the intersection numbers of a distance-regular graph.

\begin{thm}[{\cite[Main Theorem]{Weng}}]
 \label{Thm:Weng}
 Suppose that $\Gamma$ is a distance-regular graph with classical parameters $(d,b,\alpha,\beta)$ and $d \geq 4$.
 Suppose that $b < -1$, $a_1 \neq 0$, and $c_2 > 1$.
 Then $b = -q$ for some prime power $q$, and one of the following options holds:
 \begin{enumerate}
  \item $\Gamma$ is the dual polar graph of the Hermitian polar space $\mh(2d-1,q^2)$, or
  \item $\Gamma$ is the Hermitian forms graph in $\FF_{q^2}^{d \times d}$, or
  \item $q$ is odd, $\alpha = - \frac{q+1}{2}$, and $\beta = -\frac{1+(-q)^d}2$.
 \end{enumerate}
\end{thm}
No distance-regular graphs of type (3) in \Cref{Thm:Weng} are known for $d > 2$.
It follows from the work of Tian et al.\ \cite{Tian:Etal}, that these graphs do not exist for $d = 3$.
Combining some of the results by Weng \cite{Weng}, this implies the non-existence for higher values of $d$, yielding the following theorem.

\begin{thm}
 \label{Thm:extra Weng}
 Suppose that $\Gamma$ is a distance-regular graph with classical parameters $(d,b,\alpha,\beta)$ and $d \geq 4$.
 Suppose that $b < -1$, $a_1 \neq 0$, and $c_2 > 1$.
 Then $b = -q$ for some prime power $q$, and one of the following options holds:
 \begin{enumerate}
  \item $\Gamma$ is the dual polar graph of the Hermitian polar space $\mh(2d-1,q^2)$, or
  \item $\Gamma$ is the Hermitian forms graph in $\FF_{q^2}^{d \times d}$.
 \end{enumerate}
\end{thm}

By a result of Vanhove \cite{Vanhove}, this also proves that certain hemisystems in Hermitian polar spaces cannot exist.

\begin{crl}
 \label{Crl:Hemi}
 The finite Hermitian polar space $\mh(2d-1,q^2)$ with $d > 2$ and $q \neq 3$ does not have hemisystems with respect to $(d-2)$-spaces.
\end{crl}

The proofs will be discussed in more detail in \Cref{Sec:Proof}.

\section{Preliminaries}
 \label{Sec:Prel}

\subsection{Hermitian forms and polar spaces}

Consider a prime power $q$, and the finite field $\FF_{q^2}$ of square order.
The map $x \mapsto \overline x = x^q$ defines a field involution of $\FF_{q^2}$.
For any matrix $M$ over $\FF_{q^2}$, let $\overline M$ denote the matrix obtained by raising each entry of $M$ to the power $q$, and let $M^\dagger$ denote $\overline M^\top$.
The matrix $M$ is called \emph{Hermitian} if $M^\dagger = M$.

For any subspace $U \leq \FF_{q^2}^n$, define
\[
 U^\perp = \sett{y \in \FF_{q^2}^n}{(\forall x \in U)(x^\dagger y = 0}.
\]
Then the map $\perp: U \mapsto U^\perp$ is an involution on the set of subspaces of $\FF_{q^2}^n$ that reverses inclusion.
We call $\perp$ a \emph{Hermitian polarity}.
A subspace $U$ is called \emph{totally isotropic} (with respect to $\perp$) if $U \leq U^\perp$.
The Hermitian polar space $\mh(n-1,q^2)$ is the geometry formed by the totally isotropic subspaces.
We will henceforth use terminology from projective geometry.
If a subspace $U$ has bases with $t$ elements, we say that $U$ has \emph{rank $t$} and \emph{(projective) dimension $t-1$}.
Subspaces of dimension 0 are called \emph{points}.
The maximum dimension of a totally isotropic subspace of $\FF_{q^2}^n$ is given by $d = \floor{n/2}-1$, and such subspaces are called \emph{generators}.
We say that $\mh(n-1,q^2)$ is a \emph{rank $(d+1)$-geometry}.
A set $\ms$ of generators of $\mh(n-1,q^2)$ is called a \emph{hemisystem with respect to the $k$-spaces} (with $0 \leq k \leq d-2$) if for every $k$-dimensional totally isotropic subspace $U$, exactly half of the generators incident with $U$ are contained in $\ms$.
It is easy to check that a hemisystem with respect to the $k$-spaces is a hemisystem with respect to the $k'$-spaces for all $k' \leq k$.
Since the number of generators is equivalent to 1 modulo $q$, a hemisystem can only exist when $q$ is odd.

\subsection{Distance-regular graphs}

Consider a graph $\Gamma$ of diameter $d$.
For any vertex $v$ of $\Gamma$, let $\Gamma_i(v)$ denote the set of vertices at distance exactly $i$ from $\Gamma$.
We call $\Gamma$ \emph{distance-regular} if there exist numbers $a_i, b_i, c_i$ for $i=0,\dots,d$, such that for every vertex $v \in \Gamma$, and every vertex $w \in \Gamma_i(v)$,
\begin{align*}
 |\Gamma_1(w) \cap \Gamma_{i}(v)| = a_i, &&
 |\Gamma_1(w) \cap \Gamma_{i+1}(v)| = b_i, &&
 |\Gamma_1(w) \cap \Gamma_{i-1}(v)| = c_i.
\end{align*}
We remark that $\Gamma$ must be regular of degree $b_0$, and that for each $i$, $a_i + b_i + c_i = b_0$.
The sequence of numbers $\{b_0, \dots, b_{d-1}; c_1, \dots, c_d\}$ is called the \emph{intersection array} of $\Gamma$.
For more on distance-regular graphs, we refer the reader to the seminal monograph by Brouwer, Cohen, and Neumaier \cite{BCN}.
In their work, the authors noted that the intersection array of the most well-known families of distance-regular graphs (including Johnsons graphs; Grassmann graphs; dual polar graphs; bilinear, alternating, and Hermitian forms graphs; and Hamming graphs) are determined by only four numbers $(d,b,\alpha,\beta)$.
More precisely, for an integer $i$ and a real number $b$, define $[i]_b = 1+b+\ldots+b^{i-1}$.
Note that $[i]_b = \frac{b^i - 1}{b - 1}$ if $b \neq 1$.
If $\Gamma$ is a distance-regular graph of diameter $d$, and its intersection array is given by
\begin{align*}
 a_i = [i]_b \left( \beta - 1 + \alpha([d]_b - [i]_b - [i-1]_b ) \right), &&
 b_i = ([d]_b - [i]_b)(\beta-\alpha [i]_b), &&
 c_i = [i]_b(1+\alpha [i-1]_b),
\end{align*}
then $\Gamma$ is said to be \emph{classical}, and have \emph{classical parameters} $(d,b,\alpha,\beta)$.
Moreover, $b$ has to be an integer with $b \neq 0, -1$.
If $b < -1$, $\Gamma$ is said to have \emph{negative type}.

\begin{ex}
 \begin{enumerate}
  \item The \emph{dual polar graph} of $\mh(2d-1,q^2)$ is the graph whose vertices are the generators of $\mh(2d-1,q^2)$, and where $\pi_1$ and $\pi_2$ are adjacent if and only if $\dim (\pi_1 \cap \pi_2) = d-2$.
  It is a distance-regular graph with two sets of classical parameters: $\left( d,-q, - q\frac{q+1}{q-1}, - q\frac{(-q)^d+1}{q-1} \right)$ and $(d,q^2,0,q)$.
  \item The \emph{Hermitian forms graph} in $\FF_{q^2}^{d \times d}$ is the graph with as vertices the Hermitian matrices of $\FF_{q^2}^{d \times d}$, where $M$ and $N$ are adjacent if and only if $\rk(M-N) = 1$.
  It is a distance-regular graph with classical parameters $\left( d, -q, -q-1, -(-q)^d-1 \right)$.
 \end{enumerate}
\end{ex}

\begin{prop}[{\cite[Corollary 6.2.2]{BCN}, \cite{Ivanov}}]
 \label{Prop:Unique parameters}
 Suppose that $\Gamma$ is a classical distance-regular graph.
 Then either it is isomorphic to the dual polar graph of $\mh(2d-1,q^2)$, in which case it has two sets of classical parameters, or it has a unique set of classical parameters.
\end{prop}

\section{Proofs of the claims}
 \label{Sec:Proof}

In order to prove \Cref{Thm:extra Weng}, it suffices to prove that Case (3) of \Cref{Thm:Weng} cannot occur.
By combining several results of Weng \cite{Weng}, we can reduce the problem of non-existence to the case $d=3$.

\begin{thm}[{\cite[Theorem 5.7, Lemma 4.10]{Weng}}]
 \label{Thm:Reduce}
 If $\Gamma$ is a distance-regular graph with classical parameters $(d,b,\alpha,\beta)$ where $d\geq 3$, $b < -1$, $a_1 \neq 0$, and $c_2 > 1$, then for every pair of vertices $v,w$ of $\Gamma$, there exists an induced subgraph $\Gamma'$ of $\Gamma$ containing $v$ and $w$ which is distance-regular with classical parameters
 \[
  \big( t,b,\alpha,\beta + \alpha([d]_b - [t]_b) \big),
 \]
 where $t$ is the distance between $v$ and $w$.
\end{thm}

In particular, if $\Gamma$ has parameters $\left(d,b=-q,\alpha=-\frac{q+1}2, \beta = - \frac{(-q)^d+1}2\right)$, then we have
\[
 \beta + \alpha ([d]_b - [t]_b)
 = - \frac{(-q)^d+1}2 - \frac{q+1}2 \left( \frac{1-(-q)^d}{q+1} - \frac{1-(-q)^t}{q+1} \right)
 = - \frac{(-q)^t+1}2.
\]
We combine this with the recent result by Tian et al.\ \cite{Tian:Etal}.
\begin{thm}[{\cite[Theorem 1.1 (4)]{Tian:Etal}}]
 \label{Thm:Tian}
 If $\Gamma$ is a distance-regular graph with classical parameters $(3,b,\alpha,\beta)$ with $b < -1$, $a_1 \neq 0$, and $c_2 > 1$, then one of the following holds:
 \begin{enumerate}
  \item $\Gamma$ is the dual polar graph of $\mh(3,q^2)$ with parameters $\left(3,-q,-q\frac{q+1}{q-1}, q \frac{q^3-1}{q-1} \right)$, or
  \item $\Gamma$ is the extended ternary Golay code graph with parameters $(3,-2,-3,8)$, or
  \item $\Gamma$ is the large Witt graph with parameters $(3,-2,-4,10)$, or
  \item $\Gamma$ is the collinearity graph of a maximal near hexagon with classical parameters $\Big(3,-a_1-1,-\frac{c_2}{a_1}-1,(a_1+1)(c_2+a_1+1)\Big)$.
 \end{enumerate}
\end{thm}

We are now ready to prove \Cref{Thm:extra Weng}.

\begin{proof}[Proof of \Cref{Thm:extra Weng}]
 Using \Cref{Thm:Weng}, it suffices to prove that a distance-regular graph $\Gamma$ with classical parameters $\Big(d,-q,-\frac{q+1}2, - \frac{1+(-q)^d}2 \Big)$, where $q$ is an odd prime power, $d \geq 4$, $a_1 \neq 0$, and $c_2 > 1$, cannot exist.
 Hence, suppose that $\Gamma$ exists.
 Note that $a_1 = \frac{q-3}2$ and $c_2 = \frac{(q-1)^2}2$ are independent of $d$.
 Then by \Cref{Thm:Reduce}, there exists a subgraph $\Gamma'$ of $\Gamma$ which is distance-regular with classical parameters $\Big(3,-q,-\frac{q+1}2,-\frac{1-q^3}2 \Big)$.
 Note that $a_1$ and $c_2$ are the same for both graphs.
 Hence, using \Cref{Prop:Unique parameters}, the set of parameters $\Big(3,-q,-\frac{q+1}2,-\frac{1-q^3}2 \Big)$ must appear in the list of \Cref{Thm:Tian}.
 Since $q$ must be odd, cases (2) and (3) are excluded.
 If case (1) occurs, then $-\frac{q+1}2 = -q \frac{q+1}{q-1}$, which implies that $q=-1$, contradicting that $q$ is the power of a prime.
 If case (4) occurs, then $-q = -a_1-1 = - \frac{q-3}2 - 1$, which again implies that $q=-1$, a contradiction.
\end{proof}

Vanhove \cite{Vanhove} showed that hemisystems with respect to $(d-2)$-spaces of $\mh(2d-1,q^2)$ give rise to classical distance-regular graphs with parameters as in \Cref{Thm:Weng} (3).

\begin{thm}[{\cite[Theorem 4]{Vanhove}}]
 \label{Thm:Vanhove}
 Suppose that $\ms$ is a hemisystem in $\mh(2d-1,q^2)$ with respect to the $(d-2)$-spaces.
 The induced subgraph of the dual polar graph of $\mh(2d-1,q^2)$ on the set of vertices $\ms$ is distance-regular with classical parameters
 \[
  \left(d, - q, - \frac{q+1}{2}, - \frac{(-q)^d+1}{2} \right).
 \]
\end{thm}

We can now easily derive \Cref{Crl:Hemi}.

\begin{proof}[Proof of \Cref{Crl:Hemi}]
 Suppose that $d >2$, $q \neq 3$, and that $\mh(2d-1,q^2)$ has a hemisystem with respect to $(d-2)$-spaces.
 Then $q \geq 5$, since $q$ needs to be odd.
 By \Cref{Thm:Vanhove}, there exists a distance-regular graph $\Gamma$ with classical parameters $\Big(d, - q, - \frac{q+1}{2}, - \frac{(-q)^d+1}{2} \Big)$.
 As observed in the proof of \Cref{Thm:extra Weng}, $a_1 = \frac{q-3}2$ and $c_2 = \frac{(q-1)^2}2$.
 Since we suppose that $q \geq 5$, this means that $a_1 \neq 0$ and $c_2 > 1$.
 Hence, $\Gamma$ cannot exist by \Cref{Thm:extra Weng}.
\end{proof}

This leaves us with some interesting questions.

\paragraph{Open Problems.}
\begin{enumerate}
 \item Hemisystems with respect to points in $\mh(2d-1,q^2)$ do exist \cite{Bayens}.
 Can one determine the maximal $k$ such that there exist hemisystems with respect to $k$-spaces in $\mh(2d-1,q^2)$ for $d \geq 3$?
 \item This paper did not deal with the case $q=3$. Does $\mh(2d-1,3^2)$ admit hemisystems with respect to the $(d-2)$-spaces?
 More generally, do distance-regular graphs with classical parameters $\Big(d,-3,-2,-\frac{(-3)^d+1}2 \Big)$ exist for $d > 2$?
 It is known that for $d=2$, there is a unique graph with these parameters, called the (Sims-)Gewirtz graph \cite{Gerwitz}.
\end{enumerate}







\section*{Acknowledgments}

We thank Paul Terwilliger for bringing the original question to our attention.

Sam Adriaensen is supported by a postdoctoral fellowship 12A3Y25N from the Research Foundation Flanders (FWO).
Sam Mattheus is supported by a postdoctoral fellowship 1267923N from the Research Foundation Flanders (FWO). The research in this note was done during a short research retreat supported
by the Scientific Research Networks of the Research Foundation Flanders (FWO) {\em Graphs, Association schemes and Geometries: structures, algorithms and computation}, and {\em Finite Geometry, Coding Theory and Cryptography}.
This research was carried out as part of the Scientific Research Network W003324N of Research Foundation Flanders (FWO).

\newcommand{\etalchar}[1]{$^{#1}$}

\end{document}